\documentclass[11pt]{article}
\usepackage{epsfig}
\usepackage{amssymb}
\usepackage{amsmath}
\usepackage{mathrsfs}
\newcommand{\be}{\begin{equation}}
\newcommand{\ee}{\end{equation}}
\newcommand{\bea}{\begin{eqnarray}}
\newcommand{\eea}{\end{eqnarray}}

\usepackage{bm}
\usepackage{color}
\usepackage{graphicx}
\usepackage{cite}

\begin{document}


\title{On Cayley algorithm for double partition}
\author{Boris Y. Rubinstein\\
Stowers Institute for Medical Research
\\1000 50th St., Kansas City, MO 64110, U.S.A.}
\date{\today}

\maketitle
\begin{abstract}
A double partition problem asks for a number of nonnegative
integer solutions to a system of two linear Diophantine
equations with integer coefficients.
Artur Cayley suggested \cite{Cayley1860} a reduction 
of a double partition to a sum of scalar
partitions with an algorithm subject to 
a set of conditions.
We show that when these conditions are not satisfied and 
the original algorithm fails
its modification solves the reduction problem.
\end{abstract}




\section{Scalar and vector restricted partition functions}
\label{intro}

\subsection{Scalar partitions}
\label{Scalar}
The problem of integer partition into a set of integers 
is equivalent to a problem of number
of nonnegative integer solutions of the Diophantine equation
\be
\sum_{i=1}^m x_i d_i = {\bf x}\cdot{\bf d} = s.
\label{coin1}
\ee
A scalar
partition function $W(s,{\bf d}) \equiv W(s,\{d_1,d_2,\ldots,d_m\})$
solving the above problem is a
number of partitions of an integer $s$ into positive integers
$\{d_1,d_2,\ldots,d_m\}$. 
The generating function for $W(s,{\bf d})$
has a form
\be
G(t,{\bf d})=\prod_{i=1}^m\frac{1}{1-t^{d_{i}}}
 =\sum_{s=0}^{\infty} W(s,{\bf d})\;t^s\;,
\label{WGF}
\ee
Introducing notation $C[t^s](f(t))$ for a coefficient of $t^s$ in the 
expansion of a function $f(t)$ we have
\be
W(s,{\bf d}) =
C[t^s]\left( \prod_{i=1}^m(1-t^{d_{i}})^{-1} \right).
\label{const1}
\ee
Sylvester proved \cite{Sylv2} a statement about splitting of the scalar partition
function (SPF) into periodic and non-periodic parts and showed that it may be presented as a sum of "waves"
\be
W(s,{\bf d}) = \sum_{j=1} W_j(s,{\bf d})\;,
\label{SylvWavesExpand}
\ee
where summation runs over all distinct factors
of the elements of the generator vector ${\bf d}$.
The wave $W_j(s,{\bf d})$ is a quasipolynomial in $s$
closely related to prime roots $\rho_j$ of unity, namely,
is a coefficient of ${t}^{-1}$ in the series expansion 
in ascending powers of $t$ of a function
\be
F_j(s,t)=
\sum_{\rho_j} \frac{\rho_j^{-s} e^{st}}{\prod_{k=1}^{m}
\left(1-\rho_j^{d_k} e^{-d_k t}\right)}\;.
\label{generatorWj}
\ee
The summation is made over all prime roots of unity
$\rho_j=\exp(2\pi i n/j)$ for $n$ relatively prime to $j$
(including unity) and smaller than $j$.
It was shown \cite{Rub04} that it is possible to
express the wave $W_j(s,{\bf d})$ as a finite sum of the Bernoulli polynomials of
higher order.

\subsection{Vector partitions}
\label{Vector}
Consider a function $W({\bf s},{\bf D})$ counting the number of integer
nonnegative
solutions ${\bf x} \ge 0$
to a linear system ${\bf D} \cdot {\bf x} = {\bf s}$, where
${\bf D}$ is a nonnegative integer $l \times m$ generator matrix.
The function $W({\bf s},{\bf D})$ is called {\it vector partition function} (VPF) as
a natural generalization of SPF to the
vector argument.

The generating function for the VPF reads
\be
G({\bf t},{\bf D})=\prod_{i=1}^m \frac{1}{1-{\bf t}^{{\bf c}_i}} =
\sum_{{\bf s}} W({\bf s},{\bf D}) {\bf t^s},
\quad
{\bf t^s} = \prod_{k=1}^l t_k^{s_k},
\quad
{\bf t}^{{\bf c}_i} = \prod_{k=1}^l t_k^{c_{ik}},
\label{WvectGF}
\ee
where ${\bf c}_i =\{c_{i1},c_{i2},\ldots,c_{il}\},\ (1 \le i \le m)$
denotes the $i$-th column of the matrix ${\bf D}= \{{\bf c}_1,{\bf c}_2,\ldots,{\bf c}_m\}$.
Note that some elements $c_{ik}$ might equal zero.
Generalizing the coefficient notation (\ref{const1})
to the case of function of several variables
find for $W({\bf s},{\bf D})$
\be
W({\bf s},{\bf D}) 
= C[{\bf t}^{\bf s}]
\left(
\prod_{i=1}^m \frac{1}{1-{\bf t}^{{\bf c}_i}}
\right).
\label{GFvect0}
\ee

Several approaches were suggested for VPF computation including 
method of residues \cite{Beck2004,Szenes2003} and geometric decomposition
into so called {\em chambers} \cite{Sturmfels1995} -- 
regions in $m$-dimensional space where the partition function
is characterized by a specific expression.
A method of VPF computation was suggested \cite{Rub06} being
a direct generalization of the approach developed in \cite{Rub04}.
To this end {\em vector} Bernoulli and Eulerian polynomials
of higher order were introduced to find explicit
expression for $W({\bf s},{\bf D})$. A drawback of this
approach is that it does not provide any mechanism to define
chamber boundaries.
As SPF computation requires just well-known functions  \cite{Rub04}
it is promising to obtain a reduction method expressing vector partition
through scalar ones.

\subsection{Sylvester-Cayley method of vector partition reduction}
\label{Sylvester-Cayley}

The problem of scalar and vector integer partitions has a long history and J.J. Sylvester 
made a significant contribution to its solution.
In addition to the splitting algorithm for SPF \cite{Sylv2}
he suggested \cite{Sylv1} an idea to reduce VPF
into a sum of scalar partitions. 
The reduction is an iterative process based on the 
variable elimination in the generating function (\ref{WvectGF}).
Sylvester considered a specific double partition problem as 
an illustration of his method and determined regions 
of a plane $\{s_1,s_2\}$ each having a unique
expression for VPF valid in this region only.
He showed that the 
expressions in the adjacent chambers coincide at their
common boundary (see also \cite{Sturmfels1995}).

This approach was successfully applied by 
A. Cayley \cite{Cayley1860} to double partitions subject
to some restrictions on the elements of matrix ${\bf D}$ --
the vectors ${\bf c}_i$ are noncollinear,
the elements of every column ${\bf c}_i$ are relatively prime
and for all elements $c_{i2}$ of the second row the inequality $c_{i2} < s_2+2, 1\le i \le m$ holds. 
It should be noted that Cayley mentions that 
the elements of the matrix ${\bf D}$ as well as $s_i$ ``being all positive integer numbers,
not excluding zero'' \cite{Cayley1860}. Direct computation shows that when the number of 
columns containing zero is larger than two, with zeros appearing in both rows of the matrix ${\bf D}$,
and also nonzero elements in such columns are larger than unity,
the Cayley method fails. 
It also fails when the column elements have the greatest common divisor (GCD) larger than unity.
These deficiencies call for a search of alternative approach of double partition reduction
to scalar partitions. 

\subsection{Partial fraction expansion algorithm}
\label{PFE}

Computation of vector partition $W({\bf s},{\bf D})$ in (\ref{GFvect0}) 
can be performed by iterative elimination
\cite {Beck2004} of $(l-1)$ variables $t_k, (2 \le k \le l)$.
Each elimination step includes partial fraction expansion (PFE)
w.r.t. the eliminated variable
with subsequent coefficient evaluation.
This step is equivalent to 
elimination of $k$-th row of augmented $l\times(m+1)$ matrix
${\bf E}= \{{\bf s}={\bf c}_0,{\bf c}_1,{\bf c}_2,\ldots,{\bf c}_m\}$ 
made of the generator matrix ${\bf D}$
and argument vector ${\bf s}$; 
the same time one of the columns of ${\bf D}$
is eliminated too. The number of
newly generated $(l-1)\times m$ matrices is equal to $m$.
This algorithm was employed in
\cite{Cayley1860} for a two-row positive matrix 
(see also \cite{Sylv1}).




\section{Cayley algorithm of double partition reduction}
\label{Cayley}

Consider the simplest vector partition case $l=2$ following 
the algorithm described in \cite{Cayley1860}.
Denote matrix ${\bf E}$ columns as
${\bf s}={\bf c}_0=\{r,\rho\}^T$ and ${\bf c}_i=\{b_i,\beta_i\}^T,\ (1 \le i \le m),$
where ${}^T$ stands for transposition.
Cayley specified  following conditions that
should be met in order to apply the algorithm \cite{Cayley1860}. First, all fractions 
$b_i/\beta_i$ should be unequal, in other words, the 
columns ${\bf c}_i$ must be pairwise linearly independent. 
It is also required that the elements of each column ${\bf c}_i$ are relatively 
prime $\gcd(b_i,\beta_i)=1$.
Finally, all elements of the second row should satisfy a condition $\beta_i < \rho+2$.

Assuming that all $b_i, \beta_i > 0$
perform PFE step to present
$G({\bf t},{\bf D})$ as sum of $m$ fractions (${\bf t} = \{x,y\}$)
\be
G({\bf t},{\bf D}) =
\prod_{i=1}^m (1-x^{b_i}y^{\beta_i})^{-1} =
\sum_{i=1}^m T_i(x,y),
\quad
T_i = \frac{A_i(x,y)}{1-x^{b_i}y^{\beta_i}},
\label{Cayley01}
\ee
where the functions $A_i$ are rational in $x$ and rational and integral (of degree $\beta_i-1$) in $y$.

\subsection{Partial fraction expansion}
\label{Cayley1}
Consider $A_1(x,y)$ -- we have for $y=y_0=x^{-b_1/\beta_1}$,
\be
A_1(x,y) = \prod_{i \ne 1}^m (1-x^{b_i}y^{\beta_i})^{-1},
\label{mCayley02}
\ee
that is 
\be
A_1(x,x^{-b_1/\beta_1}) = \prod_{i \ne 1}^m (1-x^{b_i-b_1 \beta_i/\beta_1})^{-1}.
\label{mCayley03}
\ee
Introduce a set of complex quantities $\omega_{1j} = \omega_{1}^j,\; 1 \le j \le \beta_1-1$, where 
$\omega_{1} = \exp(2\pi i/\beta_1)$.
Multiply both the numerator and denominator of the fraction in r.h.s. in (\ref{mCayley03})
by 
$$
S_1(x,y_0) = \prod_{i \ne 1}^m  \prod_{j = 1}^{\beta_1-1} (1-\omega_{1j}x^{b_i}y_0^{\beta_i}) = 
\prod_{i \ne 1}^m  \prod_{j = 1}^{\beta_1-1} (1-\omega_{1j}x^{b_i-b_1 \beta_i/\beta_1}).
$$
The denominator turns into 
\be
\Pi_1(x) =  \prod_{i \ne 1}^m   (1-x^{\beta_1 b_i-b_1 \beta_i}),
\label{mCayley04}
\ee
while the numerator reads
\be
S_1(x,y_0) =\Pi_1(x)A_1(x,y_0) = \sum_{k = 0}^{\beta_1-1} A_{1k}(x) x^{-k b_1/\beta_1} 
= \sum_{k = 0}^{\beta_1-1} A_{1k}(x) y_0^k,
\label{mCayley05}
\ee
where $A_{1k}$ are rational functions in $x$.
This leads to 
$$
A_1(x,y_0) =\frac{1}{\Pi_1(x)}  \sum_{k = 0}^{\beta_1-1} A_{1k}(x) y_0^k,
$$
and we obtain for $T_1(x,y)$
\be
T_1(x,y) = \frac{A_1(x,y)}{1-x^{b_1}y^{\beta_1}} = 
\frac{1}{\Pi_1(x)(1-x^{b_1}y^{\beta_1})} \sum_{k = 0}^{\beta_1-1} A_{1k}(x) y^k \;.
\label{mCayley06c}
\ee

\subsection{Contribution evaluation}
\label{Cayley2}
Our goal is to find a contribution
$
C_{r,\rho}^1 = C[x^r y^{\rho}](T_1(x,y)).
$
Cayley employs a relation \cite{Cayley1860}
\be
C[x^r y^{\rho}] \left( \frac{A_1(x,y)}{1-x^{b_1}y^{\beta_1}} \right) = 
C[x^r y^{\rho}] \left( \frac{A_1(x,y)}{1-x^{b_1/\beta_1}y} \right),
\label{mCayley07}
\ee
that allows to write 
\be
\frac{A_1(x,y)}{1-x^{b_1/\beta_1}y} = U_{\beta_1-2} + \frac{A_1(x,x^{-b_1/\beta_1})}{1-x^{b_1/\beta_1}y},
\label{mCayley08}
\ee
where $U_{\beta_1-2}$ is ``a rational and integral function of the degree $\beta_1-2$ in $y$'' \cite{Cayley1860}.
Introduce
$\delta y = y - y_0$, and
perform a sequence of transformations
\bea
S_1(x,y) &=& \sum_{k = 0}^{\beta_1-1} A_{1k}(x) y^k
= \sum_{k = 0}^{\beta_1-1} A_{1k}(x) (y_0+\delta y)^k
\nonumber  \\
&=& \sum_{k = 0}^{\beta_1-1} \sum_{m=0}^{k}  A_{1k}(x) C(k,m)  y_0^m (\delta y)^{k-m}
\nonumber  \\
&=& \sum_{k = 0}^{\beta_1-1} A_{1k}(x) y_0^k
+ \sum_{k = 1}^{\beta_1-1} \sum_{m=0}^{k-1}  A_{1k}(x) C(k,m)  y_0^m (\delta y)^{k-m}
\nonumber  \\
&=& S_1(x,y_0)
+ \sum_{k = 1}^{\beta_1-1} \sum_{m=0}^{k-1}  A_{1k}(x) C(k,m)  y_0^m (y-x^{-b_1/\beta_1})^{k-m} = 
\nonumber  \\
&=& S_1(x,x^{-b_1/\beta_1})
+  \hat U,
\label{mCayley08a}
\eea
where $C(k,m)$ denotes binomial coefficient and $\hat U$ is an integral function of the degree $\beta_1-1$ in $y$
\bea
\hat U  &=& \sum_{k = 1}^{\beta_1-1} \sum_{m=0}^{k-1}  A_{1k}(x) C(k,m)  y_0^m (y-x^{-b_1/\beta_1})^{k-m}
\nonumber \\
&=& 
\sum_{k = 1}^{\beta_1-1} \sum_{m=0}^{k-1}  (-1)^{m-k} A_{1k}(x) C(k,m) x^{-kb_1/\beta_1} (1-x^{b_1/\beta_1}y)^{k-m}.
\label{mCayley08b2}
\eea
Now return to $A_1(x,y)$ in  (\ref{mCayley07}) and find for $U_{\beta_1-2}$ in (\ref{mCayley08}) 
$$
U_{\beta_1-2} = \frac{1}{\Pi_1(x)} \cdot  \frac{\hat U}{{1-x^{b_1/\beta_1}y}} = 
\frac{1}{\Pi_1(x)} \sum_{k = 1}^{\beta_1-1} \sum_{m=0}^{k-1}  (-1)^{m-k}
A_{1k}(x) C(k,m)  x^{-kb_1/\beta_1} (1-x^{b_1/\beta_1}y)^{k-1-m}\;,
$$
so that $U_{\beta_1-2}$ is 
an integral function of the degree $\beta_1-2$ in $y$
and the assumption $\beta_i < \rho + 2$ allows to drop this term from further consideration.
Then using (\ref{mCayley08}) we have
\bea
C_{r,\rho}^1 &=&C[x^r y^{\rho}] \left( \frac{A_1(x,y)}{1-x^{b_1/\beta_1}y} \right)  =
C[x^r y^{\rho}] \left( \frac{A_1(x,x^{-b_1/\beta_1})}{1-x^{b_1/\beta_1}y} \right) 
\nonumber  \\
&=&
C[x^r] \left(x^{\rho b_1/\beta_1}  A_1(x,x^{-b_1/\beta_1})\right) =
C[x^{r-\rho b_1/\beta_1}] \left(A_1(x,x^{-b_1/\beta_1})\right)
\label{mCayley09a}  \\
&=&
C[x^{r\beta_1-\rho b_1}] \left(A_1(x^{\beta_1},x^{-b_1}) \right).
\nonumber  
\eea
From (\ref{mCayley03}) we obtain 
$$
A_1(x^{\beta_1},x^{-b_1}) = \prod_{i \ne 1}^m (1-x^{b_i \beta_1-b_1 \beta_i})^{-1},
$$
and arrive at 
\be
C_{r,\rho}^1 = C[x^{r\beta_1-\rho b_1}] 
\left( \prod_{i \ne 1}^m (1-x^{b_i \beta_1-b_1 \beta_i})^{-1} \right)  
=  C[x^{r\beta_1-\rho b_1}] \left(1/\Pi_1(x) \right) .
\label{mCayley09}
\ee

\subsection{Double partition as a sum of scalar partitions}
\label{Cayley4}
Repeating computation in Sections \ref{Cayley1} and \ref{Cayley2} for each $T_i(x,y),\ 1 \le i \le m,$ we obtain 
for double partition function using (\ref{mCayley09})
$$
W(\bm s, {\bf D}) = \sum_{k=1}^m C_{r,\rho}^k = 
\sum_{k=1}^m 
C[x^{r\beta_k-\rho b_k}] 
\left( \prod_{i \ne k}^m (1-x^{b_i \beta_k-b_k \beta_i})^{-1} \right).
$$
Recall the SPF definition (\ref{const1}) and rewrite it as 
$$
W(s, {\bf d}) = C[t^s] \left( \prod_{i = 1}^m (1-t^{d_i})^{-1} \right).
$$
It allows to obtain an expression of double partition as a sum of 
scalar partitions
\be
W({\bf s},{\bf D}) = \sum_{i=1}^m W^2_i({\bf s}) =
\sum_{i=1}^m W(L_i,{\bf d}_i),
\quad
L_i = r\beta_i - b_i\rho,
\quad
d_{ij} = b_j\beta_i-b_i\beta_j,\ j\ne i.
\label{mCayley10}
\ee
This compact expression is the main result of Cayley algorithm presented in \cite{Cayley1860}.
Introducing $2 \times 2$ matrices made of ${\bf c}_i$ and the columns of augmented matrix 
${\bf E}_i=\{{\bf c}_0,{\bf c}_1,\ldots,{\bf c}_{i-1},{\bf c}_{i+1},\ldots,{\bf c}_m\}$ 
obtained from ${\bf E}$ by removal of ${\bf c}_i$
\be
{\bf D}_{i0} = \{{\bf c}_0={\bf s},{\bf c}_i\},
\quad
{\bf D}_{ij} = \{{\bf c}_j,{\bf c}_i\},
\quad
j \ne i,
\label{mCayley010}
\ee
we observe that 
$L_i={\cal D}_{i0}$ and elements of ${\bf d}_i$ are given by
$d_{ij}={\cal D}_{ij}$, where ${\cal D}_{ij} = \det {\bf D}_{ij}$.

Introduce an operation $\mathscr{R}_i$ acting on $2\times (m+1)$ augmented matrix ${\bf E}$
as follows -- first ${\bf E}$ is split into the column ${\bf c}_i$ and the matrix ${\bf E}_i$,
and then $m$ determinants ${\cal D}_{ij}$ are computed to form 
vector argument ${\cal E}_{i}=\mathscr{R}_i[{\bf E}]$ 
of the scalar partition
\be
W_i^2({\bf s}) =W({\cal E}_{i}),
\quad
{\cal E}_{i}=\mathscr{R}_i[{\bf E}] = \{{\cal D}_{ij}\},
\quad 
0 \le j \le m,\ j \ne i.
\label{mCayley100}
\ee
Linear independence of columns ${\bf c}_i$ implies that all elements
in the generator sets ${\bf d}_i$ are nonzero, but some of these might be negative
(say, $d_{ij_k}<0$ for $1 \le k \le K$). Noting that 
$(1-t^{-a})^{-1} = -t^{a}(1-t^{a})^{-1}$ we find from (\ref{mCayley10})
\be
W(L_i,{\bf d}_i) =
(-1)^{K} W(L_i + \sum_{k=1}^{K} d_{ij_k},|{\bf d}_i|),
\quad
|{\bf d}_i| = \{|d_{ij}|\}.
\label{mCayley10b}
\ee
The solution (\ref{mCayley10}) is not unique as we can
eliminate the first row of ${\bf E}$ and
obtain 
\be
W({\bf s},{\bf D}) = 
\sum_{i=1}^m W(L'_i,{\bf d}'_i),
\quad
L'_i = -L_i,
\quad
{\bf d}' = -{\bf d}.
\label{mCayley10inv}
\ee
As the term $W(L_i,{\bf d}_i)$ in (\ref{mCayley10}) and its counterpart $W(L'_i,{\bf d}'_i)$ in (\ref{mCayley10inv})
contribute only for nonnegative $L_i$ and $L'_i$ respectively, we observe that the 
terms $W(L_i,{\bf d}_i)$ and $W(L'_i,{\bf d}'_i)$ belong 
to two adjacent chambers separated by the line $L_i=0$ on which they coincide.



\section{General case of double partition reduction}
\label{General}

It was underlined above that reduction of double partition 
into a sum of scalar partitions (\ref{mCayley10}) is possible when 
several conditions on the elements of generator matrix are met.
In this Section we present an alternative approach that allows to 
drop these conditions and obtain less compact but equivalent 
reduction to scalar partitions.

\subsection{Partial fraction expansion}
\label{General1}
Consider expansion of $S_1$ in (\ref{mCayley05}) into rational functions $A_{1k}(x)$.
To find these functions write 
\be
S_1(x,y_0) =  \prod_{i \ne 1}^m  \frac{(1-x^{\beta_1 b_i}y_0^{\beta_1\beta_i})}{(1-x^{b_i}y_0^{\beta_i})}
= \prod_{i \ne 1}^m 
 \left( 
\sum_{k_i=0}^{\beta_1-1} x^{k_i b_i}y_0^{k_i \beta_i}
\right).
\label{mCayley051}
\ee
Introduce $(m-1)$-dimensional vectors
\be
\bm b'_1 = \{b_2,b_3,\ldots,b_m\},\
\bm \beta'_1 = \{ \beta_2, \beta_3,\ldots,\beta_m\},\
\bm K'_p = \{k_2,k_3,\ldots,k_m\},
\label{mCayley33b}
\ee
with $0 \le k_i \le \beta_1-1$.
Expanding r.h.s. of (\ref{mCayley051}) we obtain a sum
\be
S_1(x,y_0) = \sum_{p=1}^{n_1} x^{\bm K'_p \cdot \bm b'_1} y_0^{\bm K'_p \cdot \bm \beta'_1},
\quad
n_1=\beta_1^{m-1},
\quad
y_0=x^{-b_1/\beta_1}.
\label{mCayley33a}
\ee
For each $\bm K'_p$ we have 
$$
x^{\bm K'_p \cdot \bm b'_1} y_0^{\bm K'_p \cdot \bm \beta'_1} = 
x^{\bm K'_p \cdot \bm b'_1 - b_1 (\bm K'_p \cdot \bm \beta'_1)/\beta_1} = x^{\nu} =
x^{j_x} y_0^{j_y},
$$
where the rational
exponent $\nu$ of $x$ is split into an integer $j_x$ and 
a fractional $(\nu - j_x)$ 
parts. The fractional part of $\nu$ 
gives an exponent $j_y$ of $y_0$,
while $x^{j_x}$ contributes to $A_{1k}(x)$ in (\ref{mCayley05}).
We observe that 
\be
j_y = (\bm K'_p \cdot \bm \beta'_1) \bmod{\beta_1}.
\label{mCayley33d}
\ee
Introduce $L_1$-norm of $m$-component vector $\vert \bm a \vert = \sum_{i = 1}^{m} a_i$ and use 
(\ref{mCayley33d}) to obtain
\be
\bm K'_p \cdot \bm \beta'_1 = \beta_1 t + j_y,
\
t = \lfloor (\bm K'_p \cdot \bm \beta'_1)/\beta_1 \rfloor,
\
0 \le t=t(j_y) \le \lfloor (B_1-j_y)/\beta_1 \rfloor,
\  B_1 = (\beta_1-1) \vert  \bm \beta'_1 \vert,
\label{mCayley33dd}
\ee 
where $\lfloor \cdot \rfloor$ denotes the greatest integer less than or equal to real number.
This leads to  
\be
x^{\bm K'_p \cdot \bm b'_1} y_0^{\bm K'_p \cdot \bm \beta'_1} = 
x^ {\bm K'_p \cdot \bm b'_1 -  b_1 \lfloor (\bm K'_p \cdot \bm \beta'_1)/\beta_1 \rfloor}
 y_0^{(\bm K'_p \cdot \bm \beta'_1) \bmod{\beta_1}}.
\label{mCayley33e}
\ee
and we have for $A_1(x,y)$
\be
A_1(x,y) = \frac{1}{\Pi_1(x)} \sum_{j_y = 0}^{\beta_1-1} A_{1,j_y}(x) y^{j_y}\;,
\label{mCayley06}
\ee
with 
\be
A_{1,j_y}(x) = \sum_{j_x=N_1^{-}}^{N_1^{+}} a_{1,j_x,j_y} x^{j_x},
\quad
N_1^{-}= \min(\bm K'_p \cdot \bm b'_1 -  b_1 \lfloor (\bm K'_p \cdot \bm \beta'_1)/\beta_1 \rfloor),
\quad
N_1^{+} = (\beta_1-1) \vert  \bm b'_1 \vert.
\label{mCayley06bb}
\ee
The integer coefficients $a_{1,j_x,j_y}$ are computed from the relation 
\be
\sum_{j_y=0}^{\beta_1-1}\sum_{j_x=N_1^{-}}^{N_1^{+}} a_{1,j_x,j_y} x^{j_x} y^{j_y} = 
\sum_{p=1}^{n_1} x^{\bm K'_p \cdot \bm b'_1 -  b_1 \lfloor (\bm K'_p \cdot \bm \beta'_1)/\beta_1 \rfloor} 
y^{(\bm K'_p \cdot \bm \beta'_1) \bmod{\beta_1}},
\quad
n_1=\beta_1^{m-1},
\label{mCayley06bbb}
\ee
and the details of the computation are presented in Appendix \ref{Coeffs2Gen}.
The relation (\ref{mCayley06}) leads to
\be
T_1(x,y) = \frac{A_1(x,y)}{1-x^{b_1}y^{\beta_1}} = 
\frac{1}{\Pi_1(x)(1-x^{b_1}y^{\beta_1})} \sum_{j_y = 0}^{\beta_1-1} A_{1,j_y}(x) y^{j_y} \;.
\label{mCayley06c}
\ee

\subsection{Contribution evaluation}
\label{Cayley3a}
Drop the assumption (\ref{mCayley07}) and find $C^1_{r,\rho}$
\be
C^1_{r,\rho} =
C[x^r]\left(\Pi_1^{-1}(x)   C[y^{\rho}]\left(S_1(x,y)/(1-x^{b_1}y^{\beta_1}) \right) \right).
\label{mCayley11}
\ee
Consider the inner term in (\ref{mCayley11}) using (\ref{mCayley05})
$$
C[y^{\rho}]\left(\frac{S_1(x,y)}{1-x^{b_1}y^{\beta_1}} \right) = 
C[y^{\rho}]  \left( \sum_{j_y = 0}^{\beta_1-1}\sum_{p=0}^{\infty}  A_{1,j_y}(x)  x^{pb_1} y^{j_y+p\beta_1}  \right)
 = \sum_{j_y = 0}^{\beta_1-1}\sum_{p=0}^{\infty}  A_{1,j_y}(x)  x^{pb_1} C[y^{\rho-j_y}]  \left( y^{p\beta_1}  \right).
$$
The condition $p\beta_1 = \rho-j_y$ on integer value of $p$ reduces the inner sum to a single term
and we find 
\be
C[y^{\rho}]\left(\frac{S_1(x,y)}{1-x^{b_1}y^{\beta_1}} \right)
=\sum_{j_y = 0}^{\beta_1-1} A_{1,j_y}(x)  x^{pb_1} 
= \sum_{j_y = 0}^{\beta_1-1} A_{1,j_y}(x)  x^{(\rho-j_y)b_1/\beta_1}
= \sum_{j_y = 0}^{\beta_1-1} A_{1,j_y}(x)  x^{\rho b_1/\beta_1} y_0^{j_y}.
\label{mCayley012}
\ee
Note that the same condition on $p$  implies that 
the sum (\ref{mCayley012}) is equivalent to a single term
\be
C[y^{\rho}]\left(S_1(x,y)/(1-x^{b_1}y^{\beta_1}) \right) = 
A_{1,j_y}(x)   x^{pb_1} =  A_{1,j_y}(x)x^{(\rho-j_y)b_1/\beta_1},
\
j_y = \rho \bmod{\beta_1},
\label{mCayley12a}
\ee
where the polynomial $A_{1,j_y}(x)$ given by
(\ref{mCayley06bb}) has integer coefficients 
$a_{1,j_x,j_y}$   
determined by (\ref{mCayley06bbb}).
It provides an expression for the contribution $C^1_{r,\rho}$ of the column ${\bf c}_1$
\be
C^1_{r,\rho} = C[x^r]\left(\Pi_1^{-1}(x) \!\!\!\sum_{j=N_1^{-}}^{N_1^{+}} \!\!a_{1,j_x,j_y} x^{j+p} \right) =\!\!
\sum_{j_x=N_1^{-}}^{N_1^{+}} \!\! a_{1,j_x,j_y}  C[x^{r-j_x-(\rho-j_y)b_1/\beta_1}]\left(\Pi^{-1}_1(x)\right).
\label{mCayley12c}
\ee
On the other hand using (\ref{mCayley012}) in (\ref{mCayley11}) we obtain
\bea
C^1_{r,\rho} &=& 
C[x^{r}]\left(\Pi_1^{-1} \sum_{j_y = 0}^{\beta_1-1} A_{1,j_y}(x) x^{\rho b_1/\beta_1} y_0^{j_y}\right)
=C[x^{r-\rho b_1/\beta_1}]\left(\Pi_1^{-1} \sum_{j_y = 0}^{\beta_1-1} A_{1,j_y}(x) y_0^{j_y}\right) 
\nonumber \\ 
&=& 
C[x^{r-\rho b_1/\beta_1}]\left(S_1(x,y_0)/\Pi_1\right) = C[x^{r-\rho b_1/\beta_1}]\left(A_1(x,y_0)\right).
\label{mCayley13a}
\eea
Recalling that $y_0 = x^{-b_1/\beta_1}$ and comparing (\ref{mCayley13a}) to (\ref{mCayley09a}) we observe
that both approaches produce the same result.
Note however that the transformation 
$$
C[x^{r-\rho b_1/\beta_1}] \left(A_1(x,x^{-b_1/\beta_1})\right)=
C[x^{r\beta_1-\rho b_1}] \left(A_1(x^{\beta_1},x^{-b_1}) \right),
$$
used in (\ref{mCayley09a}) allowing  to obtain (\ref{mCayley09}) and (\ref{mCayley10})
fails when $b_1=0, \beta_1 > 1$. 
It is important to underline that the presented algorithm does not impose any restrictions on the 
value of $\rho$. 

Now using the reasoning in Section \ref{Cayley4} we write 
\bea
&&C^1_{r,\rho} = \bar W^2_1({\bf s}) = 
\sum_{j_x=N_1^{-}}^{N_1^{+}} a_{1,j_x,j_y}
W(r-j_x-(\rho-j_y)b_1/\beta_1,{\bf d}_1),
\label{mCayley13b} \\
&&j_y = \rho \bmod{\beta_1},
\quad
N_1^{+} = (\beta_1-1) \vert  \bm b'_1 \vert,
\quad
d_{1j} = b_j\beta_1-b_1\beta_j,\ j\ne 1.
\nonumber
\eea 

The vector ${\bf d}_1$ coincides with that of in (\ref{mCayley10}) and its elements 
can be computed a $2 \times 2$ determinants ${\bf D}_{ik}$ as shown in (\ref{mCayley010}).
The only difference is in the first symbolic argument of SPFs in (\ref{mCayley10c}).
We observe that it also can be written as a determinant of
$\{({\bf s} - {\bf j})/\beta_1,{\bf c}_1\}$, where ${\bf j} = \{j_x,j_y\}^T$.
Introduce a modified augmented matrix 
${\bf E}_1({\bf j})$
and write $W^2_1$ in (\ref{mCayley13b}) as 
$$ 
\bar W_1^2({\bf s}) =\!\! \sum_{j_x=N_1^{-}}^{N_1^{+}} \!\! a_{1,{\bf j}} W(\mathscr{R}_1[{\bf E}_1({\bf j})]),
\quad
{\bf E}_1({\bf j})=\{({\bf s} - {\bf j})/\beta_1,{\bf c}_1,{\bf c}_{2},\ldots,{\bf c}_m\},
$$ 
that can be presented in more general form
\bea
&&\bar W_1^2({\bf s}) = 
\!\! \sum_{{\bf j}={\bf N}_1^{-}}^{{\bf N}_1^{+}}\!\! a_{1,{\bf j}} W(\mathscr{R}_1[{\bf E}_1({\bf j})])
\delta_{j_y, \rho \bmod{\beta_1}},
\quad
{\bf j} = \{j_x,j_y\}^T,
\label{mCayley10c4} \\
&&
{\bf N}_1^{-} = \{ \min(\bm K'_p \cdot \bm b'_1 -  b_1 \lfloor (\bm K'_p \cdot \bm \beta'_1)/\beta_1 \rfloor),0\}^T,
\quad
{\bf N}_1^{+} = (\beta_1-1) \{\vert  \bm b'_1 \vert,1 \}^T.
\nonumber 
\eea
Comparing (\ref{mCayley10c4}) to (\ref{mCayley100}) we
present $\bar W_1^2$ as a weighted sum of $W_1^2$ with 
shifted argument
\be
\bar W_1^2({\bf s}) = \sum_{{\bf j}={\bf N}_1^{-}}^{{\bf N}_1^{+}} a_{1,{\bf j}} 
W_1^2(({\bf s} - {\bf j})/\beta_1)\delta_{j_y, \rho \bmod{\beta_1}},
\quad
\sum_{{\bf j}={\bf N}_1^{-}}^{{\bf N}_1^{+}} a_{1,{\bf j}} = \beta_1^{m-1}.
\label{mCayley10c5} 
\ee

\subsection{Column with nonrelatively prime elements}
\label{GCD}


The Cayley reduction method for double partition fails when at least one of the columns, say the first column
${\bf c}_1$, has GCD of its elements larger than unity $\gcd({\bf c}_1) = \gcd(b_1,\beta_1) = g_1 > 1,$
and we write $b_1=g_1 b^{\star},\ \beta_1=g_1\beta^{\star}$. 
This case can be treated as discussed above in Sections \ref{General1} and \ref{Cayley3a} leading to
$C_{r,\rho}^1$ in (\ref{mCayley12c}).
The only difference that this expansion cannot be reduced to a single term 
as in (\ref{mCayley13a}). The reason 
is that $y_0=x^{-b_1/\beta_1} =x^{-b^{\star}/\beta^{\star}}$ and thus 
upper limit in the sum in (\ref{mCayley13a}) would equal not $\beta_1$ but $\beta^{\star}$,
thus preventing a desired compactification.

\subsection{Double partition as a sum of scalar partitions}
\label{General4}
The solution (\ref{mCayley12c}) extended to other terms in PFE 
gives an expanded form of double partition equivalent to (\ref{mCayley10})
\bea
&&W({\bf s},{\bf D}) = \sum_{i=1}^m \bar W^2_i({\bf s}),
\quad
\bar W_i^2({\bf s}) = \sum_{j_x=N_i^{-}}^{N_i^{+}} a_{i,j_x,j_y}
W(r-j_x-(\rho-j_y)b_i/\beta_i,{\bf d}_i),
\label{mCayley10c} \\
&&
N_i^{+} = (\beta_i-1) \vert  \bm b'_i \vert,
\quad
d_{ij} = b_j\beta_i-b_i\beta_j,\ j \ne i,
\nonumber \\
&& 
N_i^{-}= \min(\bm K'_p \cdot \bm b'_i -  b_i \lfloor (\bm K'_p \cdot \bm \beta'_i)/\beta_i \rfloor),
\ \bm K'_p = \{k_1,k_2,\ldots,k_{i-1},k_{i+1},\ldots,k_m\},
\nonumber
\eea
where the vectors $\bm b'_i, \bm \beta'_i$
are given by
$$
\bm b'_i = \{b_1,b_2,\ldots,b_{i-1},b_{i+1},\ldots,b_m\},
\quad
\bm \beta'_i = \{\beta_1,\beta_2,\ldots,\beta_{i-1},\beta_{i+1},\ldots,\beta_m\}.
$$
Using the notation introduced in (\ref{mCayley10c4}) present $\bar W^2_i({\bf s})$ in (\ref{mCayley10c}) as 
\bea
&&\bar W_i^2({\bf s})  
= \sum_{{\bf j}={\bf N}_i^{-}}^{{\bf N}_i^{+}} a_{i,{\bf j}} 
W_i^2(({\bf s} - {\bf j})/\beta_i)
\delta_{j_y, \rho \bmod{\beta_i}},
\quad
{\bf j} = \{j_x,j_y\}^T,
\label{mCayley10cc} \\
&&
{\bf N}_i^{-} = \{ \min(\bm K'_p \cdot \bm b'_i -  b_i \lfloor (\bm K'_p \cdot \bm \beta'_i)/\beta_i \rfloor),0\}^T,
\quad
{\bf N}_i^{+} = (\beta_i-1)\{\vert  \bm b'_i \vert,1\}^T.
\nonumber 
\eea


Note that when {\it columns ${\bf c}_i$ are noncollinear} the contributions $C^i_{r,\rho}$ 
can be computed independently of each other and in this case
one can write an expression for {\it double partition as a mixture of terms $W^2_i$ and $\bar W^2_j$}.

\subsection{Double partition with collinear columns}
\label{Collinear}


The strongest condition for Cayley algorithm applicability is the
linear independence of the generator matrix columns.
It appears that in case when it fails a double partition can be
reduced to a superposition of scalar partition convolutions.

Assume that vectors corresponding to the first $n$ columns of the generator matrix {\bf D} are parallel
and rewrite the linear system ${\bf D} \cdot {\bf x} = {\bf s}$ as follows
\be
\sum_{i=1}^m {\bf c}_i x_i = {\bf s},
\quad
{\bf c}_i = u_i {\bf c}, 
\quad
{\bf c} = \{b,\beta\},
\quad
1 \le i \le n < m,
\label{noncol01}
\ee
where $u_i$ are positive integers.
The problem (\ref{noncol01}) is  equivalent to 
\be
\sum_{i=n+1}^m {\bf c}_i x_i = {\bf s}-l {\bf c} ,
\quad
\sum_{i=1}^n u_i x_i = l, 
\quad
0 \le l \le l_{max} = \min(r/b,\rho/\beta).
\label{noncol02}
\ee
Introduce a vector ${\bf u} = \{u_1,u_2,\ldots,u_n\}$ and a new matrix with $(m-n)$ columns
${\bf D}_{n+1} = \{{\bf c}_{n+1},\ldots,{\bf c}_m\},$
for which the corresponding double partition $W({\bf D}_{n+1},{\bf s}-l {\bf c})$ 
admits a reduction to SPFs either as a sum of $W_i^2$ or a mixture of $W^2_i$ and $\bar W^2_j$.
The number of solutions of (\ref{noncol02}) for given value of $l$  is equal to 
a number of solutions of the first equation given by double partition $W({\bf D}_{n+1},{\bf s}-l {\bf c})$
multiplied by a number of solutions of the second equation, {\it i.e.}, scalar partition $W(l,{\bf u})$.
Then the double partition $W({\bf D},{\bf s})$ is equivalent 
to a convolution
\be
W({\bf D},{\bf s}) = 
\sum_{l=0}^{l_{max}} W(l,{\bf u})W({\bf D}_{n+1},{\bf s}-l {\bf c}),
\quad
l_{max} = \min(r/b,\rho/\beta).
\label{noncol04}
\ee


\section{Double partitions with zero column}
\label{Zero}


Computation of contribution to double partition related to
a column ${\bf c}_1$ with zero element can be divided into three separate cases:
a) $b_1 > 0, \ \beta_1 = 0$,
b) $b_1 = 0, \ \beta_1 = 1$,
c) $b_1 = 0, \ \beta_1 > 1$.


\subsection{$b_1 > 0, \ \beta_1 = 0$}
\label{Zero_drop}

In this case application of (\ref{mCayley09})
gives 
$$
C_{r,\rho}^1 = C[x^{-\rho b_1}] \left(1/\Pi_1(x) \right) = 0,
\quad
\Pi_1(x) = \prod_{i\ne 1}^m (1-x^{-b_1 \beta_i}),
$$
as the coefficient of $x$ with negative exponent vanishes.
The corresponding scalar partition function reads
\bea
&&C_{r,\rho}^1= W(-\rho b_1,\{-b_1 \beta_2,-b_1 \beta_3,\ldots,-b_1 \beta_m\}) = 
(-1)^{m-1}W(-b_1\rho-b_1\sum_{i=2}^m \beta_i,\{b_1\beta_2,b_1\beta_3,\ldots,b_1\beta_m\})
\nonumber \\
&&=
(-1)^{m-1}W(-\rho-\sum_{i=2}^m \beta_i,\{\beta_2,\beta_3,\ldots,\beta_m\}) = 0,
\label{mCayley10dd}
\eea
where the last expression is obtained by cancelling the common factor $b_1$
in the elements of the augmented single row matrix.
Thus all columns with $\beta_i = 0$ {\it do not} contribute into 
the sum in (\ref{mCayley10}) and (\ref{mCayley10c}). 
This result is easy to understand by noticing that for every column with $\beta_i = 0$ the 
corresponding factor in (\ref{Cayley01}) does not depend on $y$ and thus 
the procedure of $y$ elimination does not involve these terms. In other words,
the corresponding term $T_i$ in the sum in (\ref{Cayley01}) is just equal to zero.

\subsection{$b_1 = 0, \ \beta_1 = 1$}
\label{Zero_beta=1}
In a particular case $\beta_1 =1$ 
we have (\ref{mCayley06})       
with $A_{10} = 1$, {\it {\it i.e.}}, $S_1(x,y)=1$ and $\Pi_1(x) = \prod_{i \ne 1}^m   (1-x^{b_i})$.
Use   (\ref{mCayley13a})    
to find
\be
C_{r,\rho}^1 = C[x^{r}]  \left(1/\Pi_1(x) \right) = 
W(r,{\bf d}_1),
\quad  
d_{1j} = b_j,\ j\ne 1.
\label{mCayley32a}
\ee
It is worth to underline that (\ref{mCayley32a}) can be also obtained by direct implementation of the original Cayley 
algorithm with $b_1=0,\; \beta_1 =1$
to find from (\ref{mCayley09}) 
\be
C_{r,\rho}^1 = C[x^{r}] 
\left( \prod_{i \ne 1}^m (1-x^{b_i})^{-1} \right)  
=  C[x^{r}] \left(1/\Pi_1(x) \right)  = W(r,{\bf d}_1),
\quad  d_{1j} = b_j,\ j\ne 1.
\label{mCayley09b}
\ee

\subsection{$b_1 = 0, \ \beta_1 > 1$}
\label{Zero_Gen}

Consider a case when a single element in the first row of matrix ${\bf D}$ equals zero and 
apply the algorithm presented above in Sections  \ref{General1}, \ref{Cayley3a}.  
Without loss of generality we can choose this {\it zero column} of ${\bf D}$
to be the first one -- ${\bf c}_1$, {\it {\it i.e.}}, $b_1 = 0$. 
In  this case (\ref{Cayley01}) turns into
\be
G({\bf t},{\bf D}) =
(1-y^{\beta_1})^{-1}
\prod_{i=2}^m (1-x^{b_i}y^{\beta_i})^{-1} = 
G_1(x,y)+\sum_{i=2}^m T_i(x,y),
\quad
G_1(x,y) = \frac{A_1(x,y)}{1-y^{\beta_1}},
\label{Cayley01b}
\ee
where $T_i(x,y)$ is given in (\ref{Cayley01})
and processed as shown above leading to the term 
$W(L_i,{\bf d}_i)$ in (\ref{mCayley10}).


Repeating the steps discussed in Sections \ref{General1}, \ref{Cayley3a}  for $b_1=0$
we obtain a particular case of (\ref{mCayley13b})

\be
C^1_{r,\rho} = \bar W^2_1 = \sum_{j_x=0}^{N_1^{+}} a_{1,j_x,j_y} 
W(r-j_x,\beta_1{\bf d}_1), 
\quad
j_y = \rho\bmod{\beta_1},
\quad
d_{1i} = b_i,\ i\ne 1.
\label{mCayley30a}
\ee
Note that scalar partition function $W(s,\beta_1{\bf d}_1) = W(s/\beta_1,{\bf d}_1)$ is nonzero only
if $s$ is divisible by $\beta_1$, and thus only terms with $j_x \equiv r\ (\bmod{\beta_1})$ 
contribute into the sum in (\ref{mCayley30a})
\be
\bar W^2_1 = \sum_{j_x=0}^{N_1^{+}} a_{1,j_x,j_y} 
W((r-j_x)/\beta_1,{\bf d}_1), \
j_y = \rho\bmod{\beta_1},\
j_x \equiv r\  (\bmod{\beta_1}),\
d_{1i} = b_i,\ i\ne 1.
\label{mCayley31a}
\ee
Computation of coefficients $a_{1,j_x,j_y}$ 
is a particular case of the general algorithm discussed 
in Appendix \ref{Coeffs2Gen}; 
in Appendix \ref{DoubleZeroExample} we present an example 
of the reduction of a double partition with zero column to a set of SPFs. 



\subsection{Multiple zero columns}
\label{Mult_2}

Generalization of the result for a single zero column discussed above 
to a case of multiple such columns with $b_i=0,\ \beta_i > 0,\ 1 \le i \le n < m,$ 
presents a particular example of linear dependent columns.
Use the results (\ref{noncol01}-\ref{noncol04}) presented in Section \ref{Collinear}
with ${\bf c} = \{0,1\},\ u_i=\beta_i,\ 1 \le i \le n < m,$ and $l_{max} = \rho$
to obtain a convolution
\be
W({\bf D},{\bf s}) = 
\sum_{l=0}^{\rho} W(l,{\bf u})W({\bf D}_{n+1},{\bf s}-l {\bf c}).
\label{Cayley84}
\ee
Note that the necessity to apply this approach arises quite rarely, namely,
when both rows of the generator matrix ${\bf D}$ have at least two zero elements.

\subsection{Alternative expression for zero column contribution}
\label{Zero_other}
Performing the elimination of the second row we obtain
\bea
W({\bf s},{\bf D}) =  \bar W^2_1 + \sum_{i=2}^{m} W^2_i
=\bar W^2_1
+\sum_{i=2}^m W(L_i,{\bf d}_i),
\quad
L_i = r\beta_i - b_i\rho,
\quad
d_{ij} = b_j\beta_i-b_i\beta_j,\ j\ne i.
\label{mCayley41}
\eea
The same time the first row elimination 
as shown in (\ref{mCayley10inv}) produces
\be
W({\bf s},{\bf D}) = 
\sum_{i=2}^m W(L'_i,{\bf d}'_i) = 
\sum_{i=2}^m W(-L_i,-{\bf d}_i) ,
\label{mCayley42}
\ee
as the term corresponding to the first column vanishes.
From (\ref{mCayley41},\ref{mCayley42}) we find an alternative expression for the
contribution of the column with $b_1=0$ as 
\be
\bar W^2_1 
= \sum_{i=2}^m \left [W(-L_i,-{\bf d}_i) - W(L_i,{\bf d}_i) \right].
\label{mCayley42a}
\ee

\section{Conclusion}
\label{Discuss}

The double partition problem subject to specific conditions
considered in \cite{Cayley1860} admits an elegant compact solution
in which the $i$-th column of the positive generator matrix  ${\bf D}$ 
leads to a single scalar partition contribution $W_i$.
This solution is obtained by elimination of a variable of
a corresponding generating function that in its turn requires
application of partial fraction expansion to the generating function.
The main steps of method are discussed in Section \ref{Cayley}; it can be applied when the generator matrix 
has  linearly independent columns, and
column elements are positive and relatively prime.

In Section \ref{General} we present  
a modification of Cayley method that instead of a single term $W_i$
produces its equivalent 
$\bar W_i$ as a weighted sum of $W_i$ with shifted argument.
Computation of the coefficients in these expressions can be 
reduced to a finite sum of double partitions with generator matrix of 
smaller size (Appendix \ref{Coeffs2Gen}).
We show that the superposition $\bar W_i$ is equivalent to the single term $W_i$
only when the specific conditions on the generator matrix elements are met
and such compactification fails when the restrictions are lifted.
An example of nonreducible superposition $\bar W_i$ is given in Section \ref{GCD} 
where we consider columns with elements that are not relatively prime.
In case of noncollinear columns each column can be processed
independently and thus a double partition can be
written as a mixture of terms $W_i$ and $\bar W_i$.
When a few columns are linear dependent a double partition leads
to a {\it convolution} of scalar partitions derived in Section \ref{Collinear}.

The case of the generator matrix  ${\bf D}$ with zero elements is 
considered in Section \ref{Zero}.
An example of a reduction of double partition with a single zero column 
is presented in Appendix \ref{DoubleZeroExample}.
Matrix with multiple zero columns discussed in Section \ref{Mult_2}
is particular case of collinear columns reducible to SPF convolution.

In conclusion we show that any double partition can be expressed through 
superposition or convolution of scalar partitions. All components of this representation
are computable using the same algorithm that makes the double partition
problem self-contained. 
As each scalar partition term $W(L_i,{\bf d}_i)$
has nonzero contribution only for $L_i \ge 0$, reduction of double partition
to SPFs allows simple determination of partition chambers.
A possibility of extension of this result to multiple partitions 
corresponding to matrices with more than two rows
remains an open question and will be discussed elsewhere.


\vskip1cm


\newpage
{\LARGE \bf Appendices}

\appendix
\section{Computation of  expansion coefficients}
\renewcommand{\theequation}{\thesection\arabic{equation}}
\setcounter{equation}{0}

\label{Coeffs2Gen}

Consider computation of integer coefficients $a_{1,j_x,j_y}$ in
(\ref{mCayley06bb}) defined through relation (\ref{mCayley06bbb})
and the vectors $\bm b', \bm \beta'$ and $\bm K'_i,\ 1\le i \le n_1,$ given by
(\ref{mCayley33b}).
For given value of $j_y$ we have from (\ref{mCayley06bbb})
$$
\sum_{j_x=N_1^{-}}^{N_1^{+}} a_{1,j_x,j_y} x^{j_x} y^{j_y} = 
\sum_{p=1}^{n_1} x^{\bm K'_p \cdot \bm b'_1 -  b_1 \lfloor (\bm K'_p \cdot \bm \beta'_1)/\beta_1 \rfloor} 
y^{(\bm K'_p \cdot \bm \beta'_1) \bmod{\beta_1}},
$$
leading to $j_y = (\bm K'_p \cdot \bm \beta'_1) \bmod{\beta_1}$, and for each value of $j_y$  we have to find 
number $a_{1,j_x,j_y}$ of vectors $\bm K'_p=\{k_2,k_3,\ldots,k_m\}$ satisfying two Diophantine equations 
with restriction on $k_i$ values
\bea
&& \bm K'_p \cdot \bm b'_1 -  b_1 \lfloor (\bm K'_p \cdot \bm \beta'_1)/\beta_1 \rfloor = j_x,
\quad
(\bm K'_p \cdot \bm \beta'_1) \bmod{\beta_1} = j_y,
\nonumber \\
&&
0 \le \bm K'_p \cdot \bm \beta'_1 \le B_1= (\beta_1-1)\vert  \bm \beta'_1 \vert,
\quad
0\le k_i \le \beta_1-1.
\label{qCayley03a}
\eea
The 
problem (\ref{qCayley03a}) thus
asks for a total number of integer solutions of $(T_1(j_y)+1)$ systems 
\be
\bm K'_p \cdot \bm b'_1 -  b_1 \lfloor(\bm K'_p \cdot \bm \beta'_1)/\beta_1 \rfloor = j_x,
\
\bm K'_p \cdot \bm \beta'_1 = j_y(t) = j_y +  \beta_1 t,
\
0 \le t \le T_1(j_y) =  \left\lfloor (B_1-j_y)/\beta_1  \right\rfloor.
\label{qCayley04a}
\ee
In the first equation in (\ref{qCayley04a}) replace $\bm K'_p \cdot \bm \beta'_1$
by $j_y + t \beta_1$ and find
$
\bm K'_p \cdot \bm b'_1 -  b_1 \lfloor (\bm K'_p \cdot \bm \beta'_1)/\beta_1 \rfloor = 
\bm K'_p \cdot \bm b'_1 - b_1 t, 
$
leading to 
\be
\bm K'_p \cdot \bm b'_1 = j_x(t) = j_x+b_1 t , 
\quad
\bm K'_p \cdot \bm \beta'_1 = j_y(t) = j_y + \beta_1  t,
\quad
0\le k_i \le \beta_1-1.
\label{qCayley05a}
\ee
Consider the problem (\ref{qCayley05a}) for a specific value of $t$. 
To include $(m-1)$ independent conditions $k_i \le \beta_1-1$ introduce a set
of $(m-1)$ additional variables $\widehat k_i,\  2 \le i \le m$ \cite{RubSylvCayl}.
Then each inequality $0 \le k_i \le \beta_1-1$ turns into an equation
$k_i + \widehat k_{i} = \beta_1-1$. 
Thus instead of two Diophantine equations in (\ref{qCayley05a}) we obtain
$(m+1)$ equations for $2(m-1)$ variables:
\be
\sum_{i=2}^m (k_i b_i + \widehat k_i \widehat b_i) = j_x(t),
\quad
\sum_{i=2}^m (k_i \beta_i + \widehat k_i \widehat \beta_i) = j_y(t), 
\quad
k_i + \widehat k_{i} = \beta_1-1,
\quad
\widehat b_i =\widehat \beta_i = 0.
\label{qCayley06a}
\ee
The system (\ref{qCayley06a}) determines the coefficient $a_{1,j_x,j_y}$ as 
a sum of vector partitions 
\be
a_{1,j_x,j_y} = \sum_{t=0}^{T_1(j_y)} W({\bf E}^{m+1}(t)),
\label{qCayley077}
\ee
where the augmented $(m+1)$-row matrix 
$${\bf E}^{m+1}(t)=\{{\bf s}^{m+1}(t),{\bf c}^{m+1}_2,{\bf c}^{m+1}_3,\ldots,{\bf c}^{m+1}_m,
\widehat{\bf c}^{m+1}_2,\widehat{\bf c}^{m+1}_3,\ldots,\widehat{\bf c}^{m+1}_m\}$$ 
is made of vectors ${\bf c}^{m+1}_i,\ \widehat{\bf c}^{m+1}_i,\ m > 2,$ and ${\bf s}^{m+1}(t),$
with vector upper index
showing its length.
The matrix ${\bf E}^{m+1}(t)$ reads
\bea
&& 
\begin{array}{ccccccccccc}
\hskip8mm{\bf s} & \hskip4.5mm{\bf c}_2 & \hskip1mm{\bf c}_3 &\ldots & {\bf c}_{m-1} & \hskip0.4mm{\bf c}_m 
&\hskip0.4mm\widehat{\bf c}_2 &\hskip0.9mm\widehat{\bf c}_3 &\hskip0.4mm\ldots  &\widehat{\bf c}_{m-1} &\hskip-2.6mm\widehat{\bf c}_m \\
\end{array}
\nonumber \\
{\bf E}^{m+1} &= &
\left(
\begin{array}{ccccccccccc}
j_x(t)  & b_2 & b_3 & \ldots & b_{m-1}  & b_m & 0 & 0 & \ldots & 0 & 0 \\
j_y(t)  & \beta_2 & \beta_3 & \ldots & \beta_{m-1} & \beta_m & 0 & 0 & \ldots & 0 & 0 \\
\beta_1 - 1 &  1 & 0 & \ldots & 0 & 0 & 1 & 0 & \ldots & 0 & 0 \\
\beta_1 - 1 &  0 & 1 & \ldots & 0 & 0 & 0 & 1 & \ldots & 0 & 0 \\
\ldots & \ldots &\ldots &\ldots &\ldots &\ldots &\ldots &\ldots &\ldots  &\ldots &\ldots \\
\beta_1 - 1 &  0 & 0 & \ldots & 1 & 0 & 0 & 0 & \ldots & 1 & 0\\
\beta_1 - 1 &  0 & 0 & \ldots & 0 & 1 & 0 & 0 & \ldots & 0 & 1
\end{array}
\right) ,
\label{qCayley06}
\eea
where we added column descriptors ${\bf s},{\bf c}_i,\widehat{\bf c}_i$ for sake of clarity.

Elimination of the last row produces nonzero contributions 
for columns ${\bf c}_m$ and $\widehat{\bf c}_m$ while all other columns have zero in the last position
and thus their contributions vanish.
Note that the last column $\widehat{\bf c}_m$ has zeros except the last unit element. 
We show in Section \ref{Zero_beta=1} that for 
zero column having last element $\beta_i=1$ its contribution to double partition
is computed using original Cayley algorithm. Similar reasoning 
applies to $\widehat{\bf c}_m$ so that we obtain
$$
W({\bf E}^{m+1}) = W(\bar{\bf E}^{m}_{0}) + W(\bar{\bf E}^{m}_{1}),
$$
where 
\bea
&& 
\begin{array}{cccccccccc}
\hskip8mm{\bf s} & \hskip4.5mm{\bf c}_2 & \hskip1mm{\bf c}_3 &\ldots & {\bf c}_{m-1} & \hskip0.4mm{\bf c}_m 
&\hskip0.4mm\widehat{\bf c}_2 &\hskip0.9mm\widehat{\bf c}_3 &\hskip0.4mm\ldots  &\widehat{\bf c}_{m-1} \\
\end{array}
\nonumber \\
\bar{\bf E}^{m}_{0} &=& 
\left(
\begin{array}{cccccccccc}
j_x(t)                 & b_2       & b_3       & \ldots & b_{m-1}      & b_m       & 0       & 0       & \ldots & 0         \\
j_y(t)            & \beta_2 & \beta_3 & \ldots & \beta_{m-1} & \beta_m & 0       & 0       & \ldots & 0         \\
\beta_1 - 1 &  1         & 0          & \ldots & 0                 & 0            & 1       & 0       & \ldots & 0        \\
\beta_1 - 1 &  0         & 1          & \ldots & 0                  & 0           & 0       & 1       & \ldots & 0        \\
\ldots         & \ldots   & \ldots    &\ldots  &\ldots             &\ldots      &\ldots & \ldots & \ldots  &\ldots  \\
\beta_1 - 1 &  0        & 0          & \ldots & 1                   & 0           & 0      & 0        & \ldots & 1            
\end{array}
\right),
\nonumber
\eea
and
\bea
&& 
\begin{array}{cccccccccc}
\hskip17mm{\bf s} & \hskip16mm{\bf c}_2 & \hskip1mm{\bf c}_3 &\ldots & {\bf c}_{m-1}  
&\hskip0.4mm\widehat{\bf c}_2 &\hskip0.9mm\widehat{\bf c}_3 &\hskip0.4mm\ldots  
&\hskip-0.6mm\widehat{\bf c}_{m-1} &\hskip-1.6mm\widehat{\bf c}_{m} \\
\end{array}
\nonumber \\
\bar{\bf E}^{m}_{1} & = & 
\left(
\begin{array}{cccccccccc}
j_x(t)  +b_m -b_m \beta_1               & b_2       & b_3       & \ldots & b_{m-1}      & 0       & 0         & \ldots    & 0  & -b_m  \\
j_y(t)  +\beta_m -\beta_m \beta_1          & \beta_2 & \beta_3 & \ldots & \beta_{m-1} & 0 & 0    & \ldots    & 0  & -\beta_m  \\
\beta_1 - 1 &  1         & 0          & \ldots & 0                 & 1            & 0       & \ldots  & 0      & 0  \\
\beta_1 - 1 &  0         & 1          & \ldots & 0                  & 0           & 1       & \ldots  & 0    & 0   \\
\ldots         & \ldots   & \ldots    &\ldots  &\ldots             &\ldots      &\ldots & \ldots & \ldots  &\ldots  \\
\beta_1 - 1 &  0        & 0          & \ldots & 1                   & 0           & 0       & \ldots  & 1    & 0  
\end{array}
\right).
\nonumber
\eea
Note that number of columns in $\bar{\bf E}^{m}_i$ is $2m-2$ -- one less as in ${\bf E}^{m+1}$. 
It is also worth to mention that the structure of the last $(m-1)$ rows of matrix ${\bf E}^{m+1}$ guarantees that value of the elements of all its rows 
except of the first two is not affected by the elimination procedure and thus to save space we would show only first two
rows of the transformed matrices.

The matrix $\bar{\bf E}^{m}_1$ has negative elements in column $\widehat{\bf c}_{m}$ and we apply 
a conversion procedure that generalizes the one employed in (\ref{mCayley10b})
$$
(1-x^{-b}y^{-\beta})^{-1} = -x^{b}y^{\beta}(1-x^{b}y^{\beta})^{-1},
$$
to obtain 
$W({\bf E}^{m+1}) =  W({\bf E}^{m}_{0}) - W({\bf E}^{m}_{1}) ,$
where ${\bf E}^{m}_0  = \bar{\bf E}^{m}_0 $ and
\bea
&& 
\begin{array}{cccccccccc}
\hskip12mm{\bf s} & \hskip10.8mm{\bf c}_2 & \hskip0.9mm{\bf c}_3 &\hskip0.4mm\ldots & {\bf c}_{m-1}  
&\hskip1mm\widehat{\bf c}_2 &\hskip1.0mm\widehat{\bf c}_3 &\hskip0.0mm\ldots  
&\hskip-1.4mm\widehat{\bf c}_{m-1} &\hskip-2.9mm\widehat{\bf c}_{m} \\
\end{array}
\nonumber \\
{\bf E}^{m}_{1} & = & 
\left(
\begin{array}{cccccccccc}
j_x(t)   -b_m \beta_1               & b_2       & b_3       & \ldots & b_{m-1}      & 0       & 0         & \ldots    & 0  & b_m  \\
j_y(t)   -\beta_m \beta_1          & \beta_2 & \beta_3 & \ldots & \beta_{m-1} & 0 & 0    & \ldots    & 0  &  \beta_m  \\
\ldots         & \ldots   & \ldots    &\ldots      &\ldots             &\ldots           &\ldots & \ldots & \ldots  &\ldots  \\
\end{array}
\right).
\nonumber
\eea
We observe that the differences between ${\bf E}^{m}_{1}$ and ${\bf E}^{m}_{0}$ include
transformation of the column ${\bf c}_m$ in ${\bf E}^{m}_{0}$ into $\widehat{\bf c}_m$ in ${\bf E}^{m}_{1}$
as well as first two elements of the argument column ${\bf s}^{m}_i$ of the matrix ${\bf E}^{m}_i$. 
Namely, ${\bf s}^{m}_{1} = {\bf s}^{m} - \beta_1 {\bf c}^{m}_m$  and ${\bf s}^{m}_{0} = {\bf s}^{m}$,
where ${\bf s}^{m}$ (${\bf c}^{m}_m$)
are obtained from ${\bf s}^{m+1}$ (${\bf c}^{m+1}_m$) belonging to matrix  ${\bf E}^{m+1}$ in (\ref{qCayley06})
by dropping the last element of the column.

Noting that both ${\bf E}^{m}_0$ and ${\bf E}^{m}_1$ have only two unit elements
in the last row
generate four matrices 
${\bf E}^{m-1}_{ij},\  i,j=0,1,$ and
find 
$$
W({\bf E}^{m}_{0}) = W({\bf E}^{m-1}_{00}) - W({\bf E}^{m-1}_{01}),
\quad
W({\bf E}^{m}_{1}) = W({\bf E}^{m-1}_{10}) - W({\bf E}^{m-1}_{11}),
$$
so that 
\be
W({\bf E}^{m+1}) = W({\bf E}^{m-1}_{00}) - W({\bf E}^{m-1}_{01}) - W({\bf E}^{m-1}_{10}) + W({\bf E}^{m-1}_{11}),
\label{qCayley07}
\ee
with 
\bea
&& 
\begin{array}{cccccccccc}
\hskip6mm{\bf s} & \hskip3.5mm{\bf c}_2 & \hskip0.9mm{\bf c}_3 &\hskip0.4mm\ldots & {\bf c}_{m-1}  & {\bf c}_{m} 
&\hskip0mm\widehat{\bf c}_2 &\hskip2.0mm\widehat{\bf c}_3 &\hskip0.0mm\ldots  
&\hskip-1.6mm\widehat{\bf c}_{m-2}  \\
\end{array}
\nonumber \\
{\bf E}^{m-1}_{00} &=& 
\left(
\begin{array}{cccccccccc}
j_x(t)                & b_2       & b_3       & \ldots & b_{m-1}      & b_m       & 0       & 0       & \ldots & 0         \\
j_y(t)            & \beta_2 & \beta_3 & \ldots & \beta_{m-1} & \beta_m & 0       & 0       & \ldots & 0         \\
\ldots         & \ldots   & \ldots    &\ldots  &\ldots             &\ldots      &\ldots & \ldots & \ldots  &\ldots  \\
\end{array}
\right),
\nonumber
\eea
\bea
&& 
\begin{array}{ccccccccc}
\hskip12mm{\bf s} & \hskip14.6mm{\bf c}_2 & \hskip0.9mm{\bf c}_3 &\hskip0.4mm\ldots & {\bf c}_{m}  
&\hskip1mm\widehat{\bf c}_2 &\hskip1.0mm\widehat{\bf c}_3 &\hskip0.0mm\ldots  
&\hskip0.0mm\widehat{\bf c}_{m-1}  \\
\end{array}
\nonumber \\
{\bf E}^{m-1}_{01} &=& 
\left(
\begin{array}{ccccccccc}
j_x(t)  -b_{m-1} \beta_1         & b_2       & b_3       & \ldots &  b_m       & 0       & 0       & \ldots & b_{m-1}         \\
j_y(t)   -\beta_{m-1} \beta_1  & \beta_2 & \beta_3 & \ldots  & \beta_m & 0       & 0       & \ldots & \beta_{m-1}         \\
\ldots         & \ldots   & \ldots    &\ldots  &\ldots             &\ldots      &\ldots & \ldots & \ldots    \\
\end{array}
\right),
\nonumber
\eea
\bea
&& 
\begin{array}{ccccccccc}
\hskip12mm{\bf s} & \hskip10.9mm{\bf c}_2 & \hskip0.9mm{\bf c}_3 &\hskip0.4mm\ldots & {\bf c}_{m-1}  
&\hskip1mm\widehat{\bf c}_2 &\hskip1.0mm\widehat{\bf c}_3 &\hskip0.0mm\ldots  
&\hskip-0.0mm\widehat{\bf c}_{m} \\
\end{array}
\nonumber \\
{\bf E}^{m-1}_{10} &=& 
\left(
\begin{array}{ccccccccc}
j_x(t)  -b_{m} \beta_1         & b_2       & b_3       & \ldots &  b_{m-1}       & 0       & 0       & \ldots & b_{m}         \\
j_y(t)   -\beta_{m} \beta_1  & \beta_2 & \beta_3 & \ldots  & \beta_{m-1} & 0       & 0       & \ldots & \beta_{m}         \\
\ldots         & \ldots   & \ldots    &\ldots  &\ldots             &\ldots      &\ldots & \ldots & \ldots   \\
\end{array}
\right),
\nonumber
\eea
and 
\bea
&& 
\begin{array}{cccccccccc}
\hskip21mm{\bf s} & \hskip18.5mm{\bf c}_2 & \hskip0.9mm{\bf c}_3 &\hskip0.2mm\ldots & {\bf c}_{m-2}  
&\hskip1mm\widehat{\bf c}_2 &\hskip1.0mm\widehat{\bf c}_3 &\hskip0.0mm\ldots  
&\hskip-0.0mm\widehat{\bf c}_{m-1} &\hskip-0.0mm\widehat{\bf c}_{m} \\
\end{array}
\nonumber \\
{\bf E}^{m-1}_{11} &=& 
\left(
\begin{array}{cccccccccc}
j_x(t)  -(b_{m-1}+b_{m}) \beta_1         & b_2       & b_3       & \ldots &  b_{m-2}       & 0       & 0       & \ldots & b_{m-1}  & b_{m}         \\
j_y(t)   -(\beta_{m-1}+\beta_{m}) \beta_1  & \beta_2 & \beta_3 & \ldots  & \beta_{m-2} & 0       & 0       & \ldots & \beta_{m-1} & \beta_{m}         \\
\ldots         & \ldots   & \ldots    &\ldots  &\ldots             &\ldots      &\ldots & \ldots & \ldots  &\ldots  \\
\end{array}
\right).
\nonumber
\eea
The argument columns ${\bf s}^{m-1}_{ij}$ of ${\bf E}^{m-1}_{ij}$ in (\ref{qCayley07}) read
$$
{\bf s}^{m-1}_{11} =  {\bf s}^{m-1} - \beta_1 ({\bf c}^{m-1}_m+{\bf c}^{m-1}_{m-1}),
\quad
{\bf s}^{m-1}_{10} =  {\bf s}^{m-1} - \beta_1 {\bf c}^{m-1}_m,
\quad
{\bf s}^{m-1}_{01} =  {\bf s}^{m-1} - \beta_1 {\bf c}^{m-1}_{m-1},
\quad
{\bf s}^{m-1}_{00} =  {\bf s}^{m-1}.
$$
Note that ${\bf s}^{m-1}$ and ${\bf c}^{m-1}_i$ represent 
first $(m-1)$ elements of the vectors ${\bf s}^{m+1}$ and ${\bf c}^{m+1}_i$ (making up 
the matrix ${\bf E}^{m+1}$) respectively.
We find after the second transformation  
\be
W({\bf E}^{m+1}) = \sum_{i,j=0}^1 (-1)^{i+j} W({\bf E}^{m-1}_{ij}),
\quad
{\bf s}^{m-1}_{ij}(t) =  {\bf s}^{m-1}(t) - \beta_1 \sum_{i,j=0}^1 (i {\bf c}^{m-1}_m+j {\bf c}^{m-1}_{m-1}).
\label{qCayley08}
\ee

It is easy to observe a pattern of iterative transformation of augmented matrices.
Introduce a vectorial index ${\bf i}(r)=\{i_1,i_2,\ldots,i_r\}$ having $r$ binary elements $i_k = 0,1$,
and denote its $L_1$-norm as $\vert {\bf i}(r)\vert = \sum_{k=1}^r i_k$.
Then we have for $r$-th transformation step
\be
W({\bf E}^{m+1}) = \!\!\! \sum_{{\bf i}(r)={\bf 0}}^{{\bf 1}} (-1)^{\vert {\bf i}(r)\vert} W({\bf E}^{m+1-r}_{{\bf i}(r)}),
\quad
{\bf s}^{m+1-r}_{{\bf i}(r)} =  {\bf s}^{m+1-r} - \beta_1 \!\!\!\! \sum_{{\bf i}(r)={\bf 0}}^{{\bf 1}} 
\sum_{k=1}^{r} i_k  {\bf c}^{m+1-r}_{m+1-k}.
\label{qCayley09}
\ee
After $(m-1)$ transformations we reduce $W({\bf E}^{m+1})$ to
\be
W({\bf E}^{m+1}) =  \!\!\!\!\!\!
\sum_{{\bf i}(m-1)={\bf 0}}^{{\bf 1}}  \!\!\!\!\!\! (-1)^{\vert {\bf i}(m-1)\vert} W({\bf E}^{2}_{{\bf i}(m-1)}),
\quad
{\bf s}^{2}_{{\bf i}(m-1)}(t)
=  {\bf s}^{2}(t) - \beta_1  \!\!\!\!\!\!\sum_{{\bf i}(m-1)={\bf 0}}^{{\bf 1}}   
\sum_{k=1}^{m-1} i_k  {\bf c}^{2}_{m+1-k},
\label{qCayley10}
\ee
with ${\bf s}^2(t) ={\bf s}^2(j_x(t),j_y(t)) = \{j_x+b_1 t,j_y+\beta_1 t\}^T$,
to computation 
of $2^{m-1}$ double partitions $W({\bf E}^{2}_{{\bf i}(m-1)})$ of augmented matrices 
\be
{\bf E}^{2}_{{\bf i}(m-1)}(j_x(t),j_y(t)) = \{{\bf s}^{2}_{{\bf i}(m-1)}(j_x(t),j_y(t)), {\bf c}^{2}_2,{\bf c}^{2}_3,\ldots,{\bf c}^{2}_m\},
\label{qCayley11}
\ee
with the columns ${\bf c}^{2}_i = \{b_i,\beta_i\}^T, \ 2 \le i \le m,$ forming a smaller generator matrix 
\be
{\bf D}^{2} = \{{\bf c}^{2}_2,{\bf c}^{2}_3,\ldots,{\bf c}^{2}_m \} = 
\left(
\begin{array}{ccccc}
 b_2       & b_3       & \ldots        & b_{m-1}        & b_m     \\
 \beta_2 & \beta_3 & \ldots         & \beta_{m-1} & \beta_m  
\end{array}
\right).
\label{qCayley12}
\ee
Introduce a reverse two-row matrix 
$$
\bar{\bf D}^{2} =\{{\bf c}^{2}_m,{\bf c}^{2}_{m-1},\ldots,{\bf c}^{2}_3,{\bf c}^{2}_2 \},
$$
define a ``scalar product'' through the discrete convolution
$$
\bar{\bf D}^{2} \cdot {\bf i}(m-1) = 
\sum_{k=1}^{m-1}  i_k  {\bf c}^{2}_{m+1-k},
$$
and rewrite ${\bf s}^{2}_{{\bf i}(m-1)}(t)$ in (\ref{qCayley10}) as
\be
{\bf s}^{2}_{{\bf i}(m-1)}(t) =  {\bf s}(t) - \beta_1  
\!\!\!\!\!\!\sum_{{\bf i}(m-1)={\bf 0}}^{{\bf 1}} \!\!\!\!\!\!\  \bar{\bf D}^{2} \cdot {\bf i}(m-1) .
\label{qCayley13}
\ee
Finally we use (\ref{qCayley10}) in (\ref{qCayley05a}) to obtain the coefficient $a_{1,j_x,j_y}$ as 
\be
a_{1,j_x,j_y} = \sum_{t=0}^{T_1(j_y)}
\sum_{{\bf i}(m-1)={\bf 0}}^{{\bf 1}} (-1)^{\vert {\bf i}(m-1)\vert} W({\bf E}^{2}_{{\bf i}(m-1)}(j_x(t),j_y(t))),
\label{qCayley14}
\ee
where $T_1(j_y)$ is defined in (\ref{qCayley04a}) and ${\bf E}^{2}_{{\bf i}(m-1)}(j_x(t),j_y(t))$ in (\ref{qCayley11}).
It is important to underline that each double partition term in (\ref{qCayley14}) is computable using 
original Cayley algorithm.

The coefficients $a_{1,j_x,j_y}$ for zero column contribution 
are computed using (\ref{qCayley14}) with $j_x(t)=j_x$.


\section{Double partition with zero column: an example}
\renewcommand{\theequation}{\thesection\arabic{equation}}
\setcounter{equation}{0}

\label{DoubleZeroExample}

Consider a double partition defined by the following augmented matrix
\be
{\bf E} = 
\left(
\begin{array}{ccccc}
s_1 & 0 & 1 & 1 & 3 \\
s_2 & 4 & 2 & 3 & 1
\end{array}
\right),
\label{zCayley00}
\ee
with a zero column ${\bf c}_1 = \{0,4\}^T$ 
and use the results presented in
Section \ref{Zero} to reduce the VPF to a sum of SPFs.
The terms $W^2_i, \ i=2,3,4$ in (\ref{mCayley41}) are computed using original Cayley 
algorithm and read
\bea
W^2_2 &=& W(2s_1-s_2-5,\{1,5,4\}),
\nonumber \\
W^2_3 &=& -W(3s_1-s_2-4,\{1,8,4\}),
\label{zCayley01} \\
W^2_4 &=& -W(s_1-3s_2-25,\{5,8,12\}).
\nonumber
\eea

The term $\bar W^2_1$ generated by elimination of zero column ${\bf c}_1$ 
requires computation of the integer coefficients $a_{1,j_x,j_y}$ in 
(\ref{mCayley30a}) where $0 \le j_y \le 3$ and $0 \le j_x \le 15$. 
We use the direct procedure based of relation (\ref{mCayley06bbb})
with $b_1=0$. 
Start with vectors $\bm b'_1, \bm \beta'_1$ and set of vectors $\bm K'_p$ 
introduced in (\ref{mCayley33b}):
$$
\bm b'_1 = \{1,1,3\},
\quad
\bm \beta'_1 = \{2,3,1\},
\quad
\bm K'_p=\{k_2,k_3,k_4\},
\quad
0 \le k_i \le 3.
$$
For every vector $\bm K'_p$ compute $x^{\bm K'_p\cdot \bm b'_1}$ and 
$j_y = (\bm K'_p \cdot \bm \beta'_1) \bmod{\beta_1}$ and for given $j_y$
collect all terms $x^{j_x}$ that allows to evaluate the coefficients $a_{1,j_x,j_y}$.
The results are presented in the table below.
$$
\begin{array}{|c|c|c|c|c|c|c|c|c|c|c|c|c|c|c|c|c|}
\hline
j_x               & 0 & 1 & 2 & 3 & 4 & 5 & 6 & 7 & 8 & 9 & 10 & 11 & 12 & 13 & 14 & 15 \\
\hline
j_y=0  & 1 & 0 & 1 & 1 & 1 & 1 & 1 & 2 & 1 & 2 &  1  &  1 &  1 &  1   &  1 & 0 \\
j_y=1  & 0 & 0 & 1 & 2 & 1 & 2 & 1 & 1 & 1 & 1 &  2 &  1  &  2 &  1   &  0 & 0 \\
j_y=2  & 0 & 1 & 1 & 1 & 1 & 1 & 2 & 1 & 2 & 1 &  1 &  1  &  1 &  1   &  0 & 1 \\
j_y=3  & 0 & 1 & 0 & 1 & 2 & 1 & 2 & 1 & 1 & 2 &  1 &  2  &  1 &  0   &  1 & 0 \\
\hline
\end{array}
$$
Using this table we write following (\ref{mCayley30a},\ref{mCayley31a})
$$
\bar W^2_1 = 
\sum_{j_x=0}^{15} a_{1,j_x,j_y} 
W(s_1-j_x,\{4,4,12\}) = 
\sum_{j_x=0}^{15} a_{1,j_x,j_y} 
W((s_1-j_x)/4,\{1,1,3\}),
$$ 
and arrive at the desired expression of the double partition
through scalar partitions
\bea
W({\bf E} ) 
&=& 
W(2s_1-s_2-5,\{1,5,4\})
-W(3s_1-s_2-4,\{1,8,4\})
\nonumber \\
&-&W(s_1-3s_2-25,\{5,8,12\})+
\sum_{j_x=0}^{15} a_{1,j_x,j_y} 
W((s_1-j_x)/4,\{1,1,3\}).
\label{zCayley02}
\eea

Employing the algorithm presented in Appendix \ref{Coeffs2Gen},
we have $m=4$, perform three iterations and generate eight two-row matrices 
${\bf E}^{2}_{{\bf i}(3)}(j_x,j_y(t))$ in (\ref{qCayley11}) with the following argument columns
\bea
&&
{\bf s}^2_{000} = \{j_x,j_y(t)\}^T,\
{\bf s}^2_{001} = \{j_x-4,j_y(t)-8\}^T,\
{\bf s}^2_{010} = \{j_x-4,j_y(t)-12\}^T,
\nonumber \\
&&
{\bf s}^2_{011} = \{j_x-8,j_y(t)-20\}^T,\
{\bf s}^2_{100} = \{j_x-12,j_y(t)-4\}^T,\
{\bf s}^2_{101} = \{j_x-16,j_y(t)-12\}^T,
\nonumber \\
&&
{\bf s}^2_{110} = \{j_x-16,j_y(t)-16\}^T,\
{\bf s}^2_{111} = \{j_x-20,j_y(t)-24\}^T,
\nonumber
\eea
while the columns 
 ${\bf c}^{2}_i = \{b_i,\beta_i\}^T, \ 2 \le i \le m,$ form a  generator matrix 
$$
{\bf D}^{2} = 
\left(
\begin{array}{ccc}
 1 & 1 & 3 \\
 2 & 3 & 1
\end{array}
\right).
$$
Applying (\ref{mCayley10}) together with (\ref{qCayley14}) we obtain 
the coefficients $a_{1,j_x,j_y}$ shown in the table above.

\end{document}